\documentclass[10pt]{amsart}

\usepackage[colorlinks=true]{hyperref}
\usepackage{upref}
\newcommand{\norm}[1]{\left\Vert#1\right\Vert}
\newcommand{\abs}[1]{\left|#1\right|}

\newcommand{\dt}{{\Delta t}}
\newcommand{\Dt}{{\Delta t}}
\newcommand{\R}{\mathbb R}
\newcommand{\Z}{\mathbb{Z}}
\newcommand{\N}{\mathbb N}
\newcommand{\Oh}{\mathcal O}
\newcommand{\dott}{\, \cdot\,}

\newcommand{\Hil}{H_{\rm per}}
\newcommand{\romnum}[1]{\MakeUppercase{\romannumeral #1}}
\newcommand{\sgn}[1]{\mathrm{sign}\left(#1\right)}
\newcommand{\lar}{\Lambda^r}

\newtheorem{theorem}{Theorem}[section]

\newtheorem{lemma}[theorem]{Lemma}

\newtheorem{remark}[theorem]{Remark}

\numberwithin{equation}{section}     
\numberwithin{figure}{section}
\numberwithin{table}{section}

\allowdisplaybreaks

\newcounter{asnr}

{\ifnum\value{asnr}=0 \stepcounter{asnr} 
  \begin{enumerate}[label=\textbf{A}.\arabic{enumi}]
    \else
    \begin{enumerate}[label=\textbf{A}.\arabic{enumi},resume] \fi}
{\end{enumerate}}

\newcounter{defnr}

{\ifnum\value{defnr}=0 \stepcounter{defnr} 
  \begin{enumerate}[label=\textbf{D}.\arabic{enumi}]
    \else
    \begin{enumerate}[label=\textbf{D}.\arabic{enumi},resume] \fi}
{\end{enumerate}}

\begin{document}
\title[Operator splitting for the BO~equation]{Operator splitting for \\  the Benjamin--Ono equation}

\author[Dutta]{R. Dutta} \address[Rajib Dutta]{\newline Department of Mathematics, University
  of Oslo, P.O.\ Box 1053, Blindern, NO--0316 Oslo, Norway}
\email[]{\href{rajib.dutta@cma.uio.no}{rajib.dutta@cma.uio.no}}

\author[Holden]{H. Holden} \address[Helge Holden]{\newline Department
  of Mathematical Sciences, Norwegian University of Science and
  Technology, NO--7491 Trondheim, Norway,\newline {\rm and} \newline
  Department of Mathematics,
  University of Oslo, P.O.\ Box 1053, Blindern, NO--0316 Oslo, Norway}
\email[]{\href{holden@math.ntnu.no}{holden@math.ntnu.no}}
\urladdr{\href{http://www.math.ntnu.no/~holden}{www.math.ntnu.no/\~{}holden}}

\author[Koley]{U. Koley} \address[Ujjwal Koley] {\newline Tata
  Institute of Fundamental Research Centre, Centre For
  Applicable Mathematics, \newline Post Bag No. 6503, GKVK Post
  Office, Sharada Nagar, Chikkabommasandra, \newline
  Bangalore 560065, India.}
\email[]{\href{ujjwal@math.tifrbng.res.in}{ujjwal@math.tifrbng.res.in}}
\urladdr{\href{http://math.tifrbng.res.in/~ujjwal/}{http://math.tifrbng.res.in/\~{}ujjwal/}}

\author[Risebro]{N. H. Risebro} \address[Nils Henrik Risebro]{\newline
   Department of Mathematics,
  University of Oslo, P.O.\ Box 1053, Blindern, NO--0316 Oslo, Norway}
\email[]{\href{nilshr@math.uio.no}{nilshr@math.uio.no}}
\urladdr{\href{http://www.mn.uio.no/math/english/people/aca/nilshr/}%
  {http://www.mn.uio.no/math/english/people/aca/nilshr/}}

\date{\today}


\thanks{Supported in part by the Research Council of Norway and the
  Alexander von Humboldt Foundation.}

\begin{abstract}
  In this paper we analyze operator splitting  for the Benjamin--Ono
  equation, $u_t= u u_x + H u_{xx}$, where $H$ denotes the Hilbert
  transform. If the initial data are sufficiently regular, we show the
  convergence of both Godunov and Strang splitting.
\end{abstract}

\maketitle


\section{Introduction}
\label{sec:intro}
In this article, we are concerned with operator splitting  for the
Benjamin--Ono equation.  The Benjamin--Ono equation models the
evolution of weakly nonlinear internal long waves. It has been derived
by Benjamin \cite{benjamin} and Ono \cite{ono} as an approximate model
for long unidirectional waves at the interface of a two-layer
system of incompressible inviscid fluids, one being infinitely
deep. In non-dimensional variables, the initial value problem
associated with the Benjamin--Ono equation reads
\begin{equation}
  \begin{cases}
    \label{eq:main}
    u_t= uu_x + H u_{xx}, \quad x\in\R, \ 0\le t \le T,& \\
    u|_{t=0}=u_0,&
  \end{cases}
\end{equation}
where $H$ denotes the Hilbert transform defined by the principle value
integral
\begin{equation*}
  H u(x) := \text{P.V.} \, \frac{1}{\pi} \int_{\R} \frac{u(x-y)}{y} \,dy.
\end{equation*}
The Benjamin--Ono equation is, at least formally, completely
integrable \cite{ablowitz} and thus possesses an infinite number of
conservation laws.  For example, the momentum and the energy, given by
\begin{align*}
  M(u):= \int u^2 \,dx, \, \, \text{and} \, \, E(u):= \frac{1}{2} \int
  \abs{D_x^{1/2} u}^2 \,dx + \frac{1}{6} \int u^3 \,dx,
\end{align*}
are conserved for solutions of \eqref{eq:main}.

We also consider the corresponding $2L$-periodic problem
\begin{equation}
  \begin{cases}
    \label{eq:main_per}
    u_t= uu_x + \Hil u_{xx}, \quad  \ 0\le t \le T,& \\
    u|_{t=0}=u_0,&
  \end{cases}
  \qquad x\in\mathbb{T},
\end{equation}
where $\mathbb{T}= \R/2 L \Z$, $u_0$ is $2L$ periodic and the periodic
Hilbert transform is defined by the principle value integral
\begin{align*}
  \Hil u(x) = \text{P.V.} \frac{1}{2L} \int_{-L}^{L}
  \cot(\frac{\pi}{2L} y) u(x-y)\,dy.
\end{align*}
The initial value problem \eqref{eq:main} has been extensively studied
in recent years. Well-posedness of \eqref{eq:main} in $H^s(\R)$ for $s
> 3$ was proved by Iorio \cite{iorio} by using purely hyperbolic
energy methods.  Then, Ponce \cite{ponce} derived a local smoothing
effect associated to the dispersive part of the equation, which
combined with compactness methods, enabled him to prove well-posedness
also for $s = 3$.

By combining a complex version of the Cole--Hopf transform with
Strichartz estimates, Tao \cite{Tao:2004} was able to show
well-posedness of the Cauchy problem \eqref{eq:main} in $H^1(\R)$.
This well-posedness was extended to $H^s(\R)$ for $s>1$ by Burq and
Planchon \cite{burq} and for $s\ge 0$ by Ionescu and Kenig
\cite{ionescu}.

In the periodic setting, Molinet \cite{molinet3} proved global well-posedness
in $H^s (\mathbb{T})$ for $s \ge 1/2$. Furthermore, he was able to improve the
global well-posedness results to $L^2(\mathbb{T})$ in \cite{Mol1}.

We employ operator splitting, i.e., the construction of an approximate
solution by concatenating the solutions of the separate problems
\begin{equation}
  \label{eq:hilbert}
  v_t=H v_{xx}
\end{equation}
and
\begin{equation}
  \label{eq:burgers}
  w_t=ww_x. 
\end{equation}
More precisely, the operator splitting method is built up as follows
\cite{Holden:book}: Consider a general partial differential equation
\begin{equation}
  \label{eq:intro2}
  u_t=C(u), \quad u|_{t=0}=u_0,
\end{equation}
where $C(u)$ is a differential operator. Furthermore, assume $C(u)$
can be written as a sum of more elementary operators, say
\begin{align*}
  C(u) = A(u) + B(u).
\end{align*}
For a positive and small time step $\Dt$ we discretize the time with
$n$ steps such that $t_n=n \Dt<T$.  Instead of solving equation
\eqref{eq:intro2} directly, we solve the two subequations
\begin{align*}
  v_t = A(v), \quad \text{and} \quad w_t = B(w),
\end{align*}
for each time step, and concatenate the solutions. The simplest form
for an operator splitting solution of \eqref{eq:intro2} is formed
solving the first subequation using the solution from the second
subequation as initial data when solving at each time step. Writing
out this procedure gives
\begin{equation}
  \begin{aligned}
    u_{n+1}=\Pi^{\Dt}(u_n)=\Phi_A^{\Dt}\circ \Phi_B^{\Dt}(u_n) =
    [\Phi_A^{\Dt}\circ \Phi_B^{\Dt}]^n (u_0),
  \end{aligned}
  \label{eq:godunov}
\end{equation}
where $u_n$ is the operator splitting solution at time $t_n$, thus
$u_n\approx u(\dott,t_n)$, and $\Phi_A^{t}(v_0)$ and $\Phi_B^{t}(w_0)$
are the exact solution operators of the above subequations at time $t$
with initial data $v_0$ and $w_0$, respectively. This is the
well-known \emph{Godunov splitting method.}

Other and more sophisticated methods for forming an operator splitting
solution of \eqref{eq:intro2} are created by solving the two
subequations for different split step sizes, and composing the
solution operators in a more complicated way. By solving one of the
subequations for half the step size composed with the solution of the
other subequation for a full time step, we obtain the \emph{Strang
  splitting method}, which is given as
\begin{equation}
  \label{eq:strang}
  u_{n+1}=\Psi^{\Dt}(u_n)=\Phi_A^{\Dt/2}\circ \Phi_B^{\Dt} \circ
  \Phi_A^{\Dt/2}(u_n) = [\Phi_A^{\Dt/2}\circ \Phi_B^{\Dt} \circ
  \Phi_A^{\Dt/2} ]^n (u_0). 
\end{equation}
For $t\in [t_n,t_{n+1})$ define $u_\Dt(t)$ by
\begin{equation*}
  u_\Dt(t)=\Pi^{t-t_n}(u_n),
\end{equation*}
in case of Godunov splitting and by
\begin{equation*}
  u_\Dt(t)=\Psi^{t-t_n}(u_n),
\end{equation*}
in case of Strang splitting.

In our case $A$ and $B$ are given by
\begin{align*}
  A(u)=H(u_{xx}), \quad \text{and} \quad B(u)=u u_x.
\end{align*}
Our main results are that the operator splitting schemes converge in $L^2$ with a rate
$\Oh(\Dt)$ for the Godunov splitting, and at a rate of $\Oh(\Dt^2)$
for the Strang splitting. However, we mention that our
method requires a well-posedness theory for the full Benjamin--Ono
equation, and cannot be used as a constructive existence theorem.  The
approach applied here has successfully been applied to a plethora of
other equations including the Korteweg--de Vries (KdV) equation, the
Schr\"odinger--Poisson, the cubic nonlinear Schr\"odinger equation,
the viscous Burgers' equation, the Benney--Lin equation, the Kawahara
equation, as well as the active scalar equation
\cite{HoldenKarlsenRisebro:1999, Lubich:2008,  HoldenKarlsenRisebroTao:2009,HoldenLubichRisebro:2011,HoldenKarlsenKarper:2011,HoldenKarlsenKarper:2013}.
 However, we stress that each equation requires its own estimates and individual treatment. In the present case both the Hilbert transform and the rather restricted well-posedness of the Benjamin--Ono equation pose new technical challenges.

The rest of the paper is organized as follows: In
Section~\ref{sec:op}, we collect well-posedness results for
\eqref{eq:main} and state results for operator splitting
schemes. Sections~\ref{sec:godunov} and~\ref{sec:strang} present the
proof of the main results for Godunov and Strang splitting,
respectively.

\section{Operator splitting}
\label{sec:op}
The upcoming analysis relies heavily on local well-posedness of
\eqref{eq:main} in $H^{s}$ in the following sense: For a given time
$T$, there exists an $R>0$ such for all $u_0\in H^s$ with $\|
u_0\|_{H^s}\le R$, there exists a unique strong solution $u\in
C([0,T], H^s)$ of \eqref{eq:main} with initial data $u_0$, and the
dependence on the initial data is locally Lipschitz continuous; i.e.,
there is a constant $K=K(R,T)<\infty$ such that for two solutions
$\widetilde u$ and $u$ corresponding to initial data $\widetilde u_0$
and $u_0$, respectively, in the $H^s$ ball of radius $R$, we have
\begin{equation}
  \label{eq:lip-flow}
  \| \widetilde u(t) - u(t) \|_{H^s} \le K \| \widetilde u_0 - u_0 \|_{H^s} 
  \quad\hbox{ for } \  0\le t \le T.
\end{equation}
Observe that this requirement says that the map taking initial data to
solution is Lipschitz continuous. Unfortunately, for the Benjamin--Ono
equation, \eqref{eq:lip-flow} is valid only for $s=0$. In fact, in
\cite{KochTzvetkov2005} it is remarked that the solution map is not
uniformly continuous from $H^s$ to $H^s$ for any $s>0$, because of the
derivative in the nonlinearity and the relatively weak smoothing
effects of the linear part of the equation.  Note, however, that the
construction in \cite{KochTzvetkov2005} does not prohibit the solution
map from being uniformly continuous or Lipschitz continuous in a
weaker topology such as $L^2$.  For the Benjamin--Ono equation we have
the following result \cite[Thm. 1.1]{Tao:2004} and
\cite[Thm. 5.3.1]{abdel}:
\begin{theorem}
  \label{theo:main}
  Let $u_0 \in H^s(\R) $ with $s > \frac{3}{2}$.  Corresponding to the
  initial data $u_0$ with $\|u_0\|_{H^s}\le R$, there exists a unique
  solution $u$ of \eqref{eq:main} with initial data $u_0$, that is,
  $u(0)=u_0$, such that
  \begin{align*}
    u \in C^k( \R_{+} ; H^{s -2k}(\R)),
  \end{align*}
  for all $k \in \N$ with $s -2k \ge -1$.  Furthermore, for another
  solution $\widetilde u(t)$ with initial data $\widetilde u_0 \in
  H^s(\R)$ such that $\|\widetilde u_0\|_{H^s}\le R$, we find
  \begin{equation*}
    \| \widetilde u(t) - u(t) \|_{L^2} \le K(R,T) \| \widetilde u_0 -
    u_0 \|_{L^2},  
    \quad\hbox{ for } \  0\le t \le T.
  \end{equation*} 
\end{theorem}
A similar result holds for the periodic case:
\cite[Thm. 1.1]{molinet3}:
\begin{theorem}
  \label{theo:main_per}
  For all $u_0 \in H^s(\mathbb{T})$ with $s \ge \frac{1}{2}$ and for
  all $T>0$, there exists a solution $u$ of the Benjamin--Ono equation
  \eqref{eq:main_per} satisfying
  \begin{align*}
    u \in C\left( [0,T]; H^s(\mathbb{T}) \right).
  \end{align*}
  Moreover, $u \in C_b(\R,L^2(\mathbb{T}))$ and the map $u_0 \mapsto
  u$ is continuous from $H^s$ into $ C\left( [0,T]; H^s(\mathbb{T})
  \right)$ and Lipschitz on every bounded set from $H_0^s$ into $
  C\left( [0,T]; H_0^s(\mathbb{T}) \right)$. Here $H_0^s(\mathbb{T})$
  denotes the closed subset of $H^s(\mathbb{T})$ with mean zero.
\end{theorem}
For Godunov splitting, we consider solutions bounded by
\begin{equation}
  \label{eq:HYP}
  \norm{u(t)}_{H^{5/2}} \le \rho <R \quad\hbox{ for } \ 0\le t \le T,
\end{equation}
in particular, $u_0\in H^{5/2}$.  We show the following result.
\begin{theorem}[First-order convergence in $L^2$]
  \label{thm:one}
  Let $u$ be the unique solution of \eqref{eq:main}, and assume that
  $u$ satisfies \eqref{eq:HYP}.  Define the Godunov approximation
  $u_\Dt$ by \eqref{eq:godunov}.  Then for any $T>0$ there is a
  $\overline{\dt}>0$ such that for $\dt\le \overline{\dt}$ and $t\le
  T$, we have
  \begin{equation*}
    \norm{u_\Dt(t) - u(t_n)}_{L^{2}} \le C_1 \dt.
  \end{equation*}
  Here, $\overline{\dt}$ and $C_1$ only depend on $\|u_0\|_{H^{5/2}}$,
  $\rho$, and $T$.
\end{theorem}
Regarding Strang splitting, we consider solutions bounded by
\begin{equation}
  \label{eq:HYP1}
  \norm{u(t)}_{H^{9/2}} \le \rho <R \quad\text{for}\quad \ 0\le t \le T.
\end{equation}
In this case, we assume that $u_0\in H^{9/2}$. Then we show the
following result.
\begin{theorem}[Second-order convergence in $L^2$]
  \label{thm:two}
  Let $u$ be the unique solution of \eqref{eq:main}, and assume that
  $u$ satisfies \eqref{eq:HYP1}.  Define the Strang approximation
  $u_\Dt$ by \eqref{eq:strang}.  Then there is a $\overline{\dt}>0$
  such that for $\dt\le \overline{\dt}$ and $t\le T$, we have
  \begin{equation*}
    \norm{u_\Dt(t) - u(t)}_{L^{2}} \le C_2 \Dt^2.
  \end{equation*}
  Here, $\overline{\dt}$ and $C_2$ only depend on $\|u_0\|_{H^{9/2}}$,
  $\rho$, and $T$.
\end{theorem}
Since the exact solution operator for Burgers' equation eventually
will produce discontinuities independently of the smoothness of the
initial data, the initial value problem for Burgers' equation is not
well posed in any Sobolev space with positive exponent. However, if
the initial values are smooth, discontinuities will not be created
instantaneously, and if you know that the solution is smooth, it is
actually smoother than you think. The precise result reads as follows.
\begin{lemma}\label{lem:burgers-0}
  Let $r\ge r_0>3/2$. Then following results hold: \\
  (i) If $\|\Phi_B^t(u_0) \|_{H^{r_0}} \le \alpha$ for $0\le t\le\dt$,
  then $\|\Phi_B^t(u_0) \|_{H^r} \le e^{C\alpha t}\| u_0 \|_{H^r}$ for
  $0\le t\le\dt$, where the constant $C$ is
  independent of $u_0$ and $\dt$. \\
  (ii) If $\| u_0 \|_{H^{r_0}} \le M$, then there exists a time
  $\bar{t}(M)>0$ such that $\| \Phi^t_B(u_0) \|_{H^{r_0}} \le C(M)$
  for $0\le t \le \bar{t}(M)$.
\end{lemma}
\begin{proof}
  For any number $r \in \R$, let $H^r(\R)$ be the Sobolev space
  consisting of all tempered distributions $f$ such that
  \begin{align*}
    \norm{f}_r = \left( \int_{\R} {\langle \xi \rangle}^{2r}
      |{\hat{f}(\xi)}|^2 \,d\xi\right)^{1/2} < \infty,
  \end{align*}
  with $\langle \xi \rangle = (1 + |\xi|^2)^{1/2}$ and $\hat{f}(\xi) =
  \mathcal{F}(f) (\xi)= \int_{\R} e^{i \xi x} f(x) \,dx$ is the Fourier
  transform of $f$.  Furthermore, we define an integral operator
  $\lar$ on tempered distributions by
  \begin{align*}
    \lar(f)=\mathcal{F}^{-1}( {\langle \xi \rangle}^r \hat{f}),
  \end{align*}
  where $\mathcal{F}^{-1}$ denotes the inverse Fourier
  transform. Since the inverse Fourier transformation preserves the
  $L^2$ norm, it is evident that $\norm{\lar(f)}_{L^2} =
  \norm{f}_{H^r}$. Moreover, it is easy to see that $\lar$ is linear,
  and commutes with the derivative, i.e., $\lar(f)_x=\lar(f_x)$.

  Let $u$ be a solution of Burgers' equation. Then
  $\lar(u_t)=\lar(uu_x)$. Taking the standard $L^2$ inner product,
  denoted $\langle\dott,\dott\rangle_{L^2}$, with $\lar u$ yields
  \begin{align*}
    \frac12 \frac{d}{dt} \norm{\lar u}_{L^2}^2&= \langle \lar u, \lar
    u_t  \rangle_{L^2}=  \langle \lar u, \lar u u_x \rangle_{L^2}\\ 
    &=\langle \lar u, u\lar u_x \rangle_{L^2} + \langle \lar u, \lar u
    u_x -u\lar u_x \rangle_{L^2}.
  \end{align*}
  The first term of the above expression can be estimated as follows
  \begin{align*}
    \abs{\langle \lar u, u\lar u_x \rangle_{L^2}}&= \Bigl|\int_{\R}
    u(\lar u)_x (\lar u) dx\Bigr|=\Bigl|\frac12 \int_{\R} u_x (\lar
    u)^2 dx\Bigr|
    \\
    &\le \norm{u_x}_{L^\infty} \norm{\lar u}^2_{L^2}\le C
    \norm{u}_{H^{r_0}}\norm{u}^2_{H^r},
  \end{align*}
  where we have used the Sobolev inequality
  \begin{align*}
    \norm{u_x}_{L^\infty} \le C \norm{u_x}_{H^{r_0-1}} \le C
    \norm{u}_{H^{r_0}},
  \end{align*}
  which holds since $r_0-1>1/2$.  The second term can be estimated by
  the Cauchy--Schwarz inequality, i.e.,
  \begin{equation*}
    \abs{\langle \lar u, \lar u u_x -u\lar u_x \rangle_{L^2}}\le \norm{\lar
      u}_{L^2} \norm{\lar u u_x -u\lar u_x}_{L^2}.
  \end{equation*}
  To proceed further, we need the following inequalites which can be
  readily verified using the mean value theorem: For $r >1$, and any
  $\xi$ and $\eta$,
  \begin{subequations}
    \begin{align}
      \abs{(1+ \xi^2)^{r/2} - (1+\eta^2)^{r/2} }&\le C \abs{\xi-\eta}
      \left[ (1+(\xi-\eta)^2)^{\frac{r-1}2}+
        (1+\eta^2)^{\frac{r-1}2}\right],\label{eq:ineq1}\\
      \abs{\eta}&\le \left(1+\eta^2\right)^{1/2}, \label{eq:ineq2}
    \end{align}
  \end{subequations}
  where $C$ is a constant. At this point, we also recall Young's
  inequality for convolutions
  \begin{align*}
    \norm{u * v}_{L^2}\le \norm{u}_{L^1} \norm{v}_{L^2}.
  \end{align*}
  With the above inequalities, we calculate
  \begin{align*}
    &\norm{\lar u u_x -u\lar u_x}_{L^2}\\
    &= \Bigl(\int_\R \Bigl( \int_\R \left((1+ \xi^2)^{r/2} -
      (1+\eta^2)^{r/2}\right) \hat{u}(\xi-\eta) \hat{u}_x(\eta)
    d\eta \Bigr)^2d\xi\Bigr)^{1/2}\\
    &\le C \Bigl(\int_\R \Bigl( \int_\R \left[
      (1+(\xi-\eta)^2)^{\frac{r-1}2}+ (1+\eta^2)^{\frac{r-1}2}\right]
    \abs{\xi-\eta} \abs{\hat{u}(\xi-\eta) \hat{u}_x(\eta) }
    d\eta \Bigr)^2d\xi\Bigr)^{1/2}\\
    &\le C\Bigl(\int_\R \Bigl( \int_\R
    \left[1+(\xi-\eta)^{\frac{r-1}{2}}\right]\abs{\xi-\eta}
    \,\abs{\hat{u}(\xi-\eta)
      \hat{u}_{x}(\eta)}\,d\eta\Bigr)^2\,d\xi\Bigr)^{1/2} \\
    &\qquad \qquad \qquad \qquad \qquad \qquad +
    C\Bigl(\int_\R\Bigl(\int_\R
    \left[1+\eta^{\frac{r-1}{2}}\right]\abs{\hat{u}_x(\xi-\eta)}\,
    \abs{\eta \hat{u}(\eta)}\,d\eta\Bigr)^2\,d\xi\Bigr)^{1/2}\\
    &\le C\Bigl(\int_\R\Bigl(\int_\R
    \left[1+(\xi-\eta)^{\frac{r}{2}}\right]\abs{\hat{u}(\xi-\eta)}\,
    \abs{\hat{u}_x(\eta)}\,d\eta\Bigr)^2\,d\xi\Bigr)^{1/2}\\
    &\qquad \qquad \qquad \qquad \qquad \qquad \quad+
    C\Bigl(\int_\R\Bigl(\int_\R
    \left[1+\eta^{\frac{r}{2}}\right]\abs{\hat{u}(\eta)}\,
    \abs{\hat{u}_x(\xi-\eta)}\,d\eta
    \Bigr)^2\,d\xi\Bigr)^{1/2}\\
    &\le
    C\norm{\hat{u}_x}_{L^1}\,\norm{(1+\xi^2)^r\hat{u}(\xi)}_{L^2}.
  \end{align*}
  For the first factor, observe that
  \begin{align*}
    \norm{\hat{u}_x}_{L^1} &= \int_{\R} \abs{\hat{u}_x (\xi)} \,d\xi
    =\int_\R \abs{\xi}\,\abs{\hat{u}(\xi)}\,d\xi\\
    &\le \Bigl( \int_\R\left(1 + \xi^2\right)^{r_0}
    \abs{\hat{u}(\xi)}^2 \,d\xi \Bigr)^{1/2} \Bigl( \int_{\R}
    \frac{\xi^2}{(1 +
      \xi^2)^{r_0}} \,d\xi \Bigr)^{1/2} \\
    &= C_{r_0} \norm{u}_{H^{r_0}}, \quad \text{for}\,\, r_0 > 3/2.
  \end{align*}
  Thus, combining the above estimates, we obtain
  \begin{equation*}
    \abs{\langle \lar u, \lar u u_x -u\lar u_x \rangle_{L^2}}\le C
    \norm{u}_{H^{r_0}}\norm{u}_{H^r}^2.
  \end{equation*}
  Therefore
  \begin{align}
    \label{eq:arg}
    \frac{d}{dt} \norm{u}_{H^r}^2\le C\norm{u}_{H^{r_0}} \norm{
      u}_{H^r}^2
  \end{align}
  which proves the first part (i) of the Lemma
  ~\ref{lem:burgers-0}. Observe that, we can also use $r_0=r$ in
  \eqref{eq:arg}, which implies
  \begin{align}
    \label{eq:arg1}
    \frac{d}{dt} \norm{ u}_{H^r}\le C \norm{ u}_{H^r}^2.
  \end{align}
  The second part (ii) of Lemma~\ref{lem:burgers-0} follows by
  comparing \eqref{eq:arg1} with the majorizing differential equation
  $y' = C y^2$.
\end{proof}

\section{Godunov splitting}
\label{sec:godunov}

In the previous subsection, we have presented several results which
now will prove useful. In what follows, we first estimate the local
error for the Godunov splitting, before we use this estimate to find a
bound for the global error.

We start by a general perturbation result. We write
$e^{tA}v=\Phi^t_A(v)$ to indicate the linearity of the flow of $A$. We
start from the variation-of-constants formula \cite[Thm. 4.2.4]{Leela}
 for $u(t)=\Phi^t_{A+B}(u_0)$,
\begin{equation}
  \label{voc}
  u(t)=e^{tA}u_0 + \int_0^t e^{(t-s)A} B(u(s))\, ds,
\end{equation}
which is just the formula $\phi(t)-\phi(0)=\int_0^t \dot\phi(s)\,ds$
for $\phi(s)=e^{(t-s)A}u(s)$.  Furthermore, we have\footnote{Here we
  introduce the second-order Taylor expansion
  $\Psi(f+g)=\Psi(f)+d\Psi(f)[g]+\int_0^1
  (1-\alpha)d^{(2)}\Psi(f+\alpha g)[g]^2d\alpha$ for an operator
  $\Psi$, see \cite[p.~29]{AmbrosettiProdi} for notation and proofs.}
\begin{equation}
  \label{voc-b}
  B(u(s)) = B(e^{sA}u_0) + \int_0^s dB(e^{(s-\sigma) A} u(\sigma))[e^{(s-\sigma) A}B(u(\sigma))]\, d\sigma,
\end{equation}
which is nothing but the formula $B(\varphi(s))-B(\varphi(0))=\int_0^s
dB(\varphi(\sigma))[\dot\varphi(\sigma)]\,d\sigma$ for
$\varphi(\sigma)=e^{(s-\sigma) A} u(\sigma)$. We insert \eqref{voc-b}
into \eqref{voc} with $t=\dt$ to obtain
\begin{equation} \label{voc1AA} u(\dt)=e^{\dt A}u_0 + \int_0^\dt
  e^{(\dt-s)A} B(e^{sA}u_0)\, ds + e_1
\end{equation}
with
\begin{equation}
  \label{e1}
  e_1 = \int_0^\dt \int_0^s e^{(\dt-s)A}dB(e^{(s-\sigma) A} u(\sigma))[e^{(s-\sigma) A}B(u(\sigma))]\, d\sigma \, ds.
\end{equation}
We next turn to results specifically for the Godunov splitting. The
main tool for proving Theorem~\ref{thm:one} is a local error estimate.
\begin{lemma}
  \label{lem:local-1}
  Assume that hypothesis \eqref{eq:HYP} holds for the solution $u(t) =
  \Phi_{A+B}^t(u_0)$ of \eqref{eq:main}.  If the initial data $u_0$ is
  in $H^{5/2}$, then the local error of the Godunov splitting
  \eqref{eq:godunov} is bounded in $L^2$ by
  \begin{equation*}
    \| \Pi^\dt(u_0)-\Phi_{A+B}^\dt(u_0) \|_{L^2} \le c_1 \dt^2,
  \end{equation*}
  where $c_1$ only depends on $\|u_0\|_{H^{5/2}}$.
\end{lemma}
\begin{proof} Set
  \begin{equation*}
    u_1 = \Pi^\dt(u_0)= e^{\dt A} \left(\Phi^\dt_B(u_0)\right).
  \end{equation*}
  The first-order Taylor expansion with integral remainder term
  \begin{equation}
    \label{taylor}
    \Phi^\dt_B (v) = v + \dt B(v) + \underbrace{\dt^2 
      \int_0^1 (1-\theta) dB(\Phi^{\theta\dt}_B
      (v))[B(\Phi^{\theta\dt}_B (v))]\, d\theta}_{e_2}
  \end{equation}
  is justified for any $v \in H^{5/2}$ and for sufficiently small
  $\dt$ by Lemma~\ref{lem:burgers-0}. We therefore obtain
  \begin{equation*}
    u_1= e^{\dt A} u_0 + \dt e^{\dt A} B(u_0) + e_2. 
  \end{equation*}
  Thus the error can be written
  \begin{equation}
    \label{err1}
    u_1-u(\dt) = \dt\, e^{\dt A} B(u_0) - \int_0^\dt e^{(\dt-s)A}
    B(e^{sA}u_0)\, ds + (e_2-e_1), 
  \end{equation}
  and therefore the principal error term is just the quadrature error
  of the rectangle rule applied to the integral over $[0,\dt]$ of the
  function
  \begin{equation}
    \label{f}
    f(s) = e^{(\dt-s)A} B(e^{sA}u_0).
  \end{equation}
  We express the quadrature error in first-order Peano form,
  \begin{equation*}
    \dt\, f(0) -\int_0^\dt f(s) \, ds = \dt^2 \int_0^1 \kappa_1(\theta)\,
    f'(\theta\dt)\, d\theta,
  \end{equation*}
  where $\kappa_1$ is the real-valued, bounded Peano kernel of the
  rectangle rule.  Thus, the $L^2$-error after one step is bounded as
  \begin{equation}
    \label{error_full}
    \norm{u_1-u(\dt)}_{L^2}  \le  (\dt)^2   \int_0^1
    \norm{\kappa_1(\theta)\, f'(\theta\dt)}_{L^2}\, d\theta +
    \norm{(e_2-e_1)}_{L^2}. 
  \end{equation}
  Next, we find that
  \begin{equation*}
    f'(s) = -e^{(\dt-s)A} [A,B](e^{sA}u_0)
  \end{equation*}
  \begin{equation} \label{eq:lie1}
    \begin{aligned}\relax
      [A,B](v) & = dA(v)[B(v)] - dB(v)[Av]  \\
      & = H[(v v_x)_{xx}] - v_x H(v_{xx}) - v H(v_{xxx})  \\
      & = H(v v_{xxx}) + 3 H(v_x v_{xx}) -v_x H(v_{xx}) - v
      H(v_{xxx}).
    \end{aligned}
  \end{equation}
  By Lemma~\ref{lem:1COM} below, we obtain the commutator bound,
  \begin{equation*}
    \norm{[A,B](v)}_{L^2} \le C \norm{v}_{H^{5/2}}^2.
  \end{equation*}
  Since $e^{tA}$ preserves the Sobolev norms, we have
  \begin{equation*}
    \langle u,Hu_{xx}\rangle_{L^2}=- \langle Hu,u_{xx}\rangle_{L^2}=-
    \langle(Hu)_{xx},u\rangle_{L^2}=- \langle Hu_{xx},u\rangle_{L^2}, 
  \end{equation*}  
  implying that $\langle u,Hu_{xx}\rangle_{L^2}=0$. It follows 
  \begin{equation*}
    \norm{f'(s)}_{L^2} \le C\norm{u_0}_{H^{5/2}}^2,
  \end{equation*}
  and hence the quadrature error is $O(\dt^2)$ in the $L^2$-norm for
  $u_0\in H^{5/2}$.  The $L^2$-norm of the remainder term $e_2-e_1$ is
  bounded by $C\dt^2$ for $u_0\in H^{5/2}$ for sufficiently small
  $\dt$ (by using Lemma~\ref{lem:burgers-0} (ii)). Specifically,
  \begin{align*}
    \norm{e_1}_{L^2}&\le \int_0^\dt \int_0^s
    \norm{e^{(\dt-s)A}dB(e^{(s-\sigma) A} u(\sigma))[e^{(s-\sigma)
        A}B(u(\sigma))]}_{L^2}\, d\sigma \, ds \\
    &\le \int_0^\dt \int_0^s\norm{\Big(\big(e^{(s-\sigma) A}
      u(\sigma)\big)\big(e^{(s-\sigma)
        A}B(u(\sigma))\big)\Big)_x}_{L^2}\, d\sigma \, ds\\
    &\le C\int_0^\dt \int_0^s\norm{u(\sigma)}_{H^{1}}
    \norm{B(u(\sigma))}_{H^{1}}\, d\sigma \, ds\\
    &\le C\int_0^\dt \int_0^s\norm{u(\sigma)}_{H^{1}}
    \norm{u(\sigma)}_{H^{1}}\ \norm{u(\sigma)}_{H^{2}}\, d\sigma \, ds\\
    & \le C\int_0^\dt \int_0^s \norm{u(\sigma)}^3_{H^{2}}\, d\sigma \,
    ds \, \le C \dt^2 R^3, \intertext{and} \norm{e_2}_{L^2}&\le \dt^2
    \int_0^1 (1-\theta)
    \norm{e^{\dt A} dB(\Phi^{\theta\dt}_B (u_0))[B(\Phi^{\theta\dt}_B
      (u_0))]}_{L^2}\, d\theta\\
    &\le \dt^2 \int_0^1 \norm{\Big((\Phi^{\theta\dt}_B
      (u_0))(B(\Phi^{\theta\dt}_B
      (u_0)))\Big)_x}_{L^2}\, d\theta\\
    &\le\dt^2 C\int_0^1 \norm{\Phi^{\theta\dt}_B (u_0)}_{H^{1}}
    \norm{B(\Phi^{\theta\dt}_B
      (u_0))}_{H^{1}}\, d\theta\\
    &\le\dt^2 C\int_0^1 \norm{\Phi^{\theta\dt}_B (u_0)}_{H^{1}}^2
    \norm{\Phi^{\theta\dt}_B
      (u_0)}_{H^{2}}\, d\theta\\
    &\le\dt^2 C\int_0^1 \norm{\Phi^{\theta\dt}_B (u_0)}^3_{H^{2}}\,
    d\theta \, \le C\dt^2 R^3.
  \end{align*}
  This completes the proof.
\end{proof}
\begin{lemma} \label{lem:1COM} If $v\in H^{5/2}(\R)$ and $[A,B](v)$ is
  given by \eqref{eq:lie1}, then
  \begin{equation}
    \norm{[A,B](v)}_{L^2} \le C\norm{v}_{H^2}\norm{v}_{H^{5/2}},
  \end{equation}
  for some constant $C$.
\end{lemma}
\begin{proof}
  From \eqref{eq:lie1} we have
  \begin{equation*} [A,B](v) = g(v)_x + 2 H(v_x v_{xx}),
  \end{equation*} 
  with
  \begin{equation*}
    g(v):= H(v v_{xx}) - v H(v_{xx}).
  \end{equation*}
  First we show
  \begin{equation*}
    \norm{g(v)_x}_2\le \norm{v_x}_2\norm{v}_{H^{5/2}}.
  \end{equation*} 
  Recall that $\mathcal{F}$ denotes the Fourier transform. Thus
  \begin{equation*}
    \mathcal{F}(g(v))(\xi)=  \romnum{1}(\xi)- \romnum{2}(\xi)
  \end{equation*}
  where
  \begin{align*}
    \romnum{1}(\xi)&:=\mathcal{F}\left(vH(v_{xx})\right)(\xi)
    =-i \int \hat{v}(\xi-\eta)\, \hat{v}(\eta)\,
    \sgn{\eta}\, |\eta|^2d\eta,\\
    \romnum{2} (\xi)&:=\mathcal{F}(H(vv_{xx}))(\xi)
    =-i\, \sgn{\xi} \int\hat{v}(\xi-\eta)\,
    \hat{v}{\eta}\, |\eta|^2d\eta.
  \end{align*}
  Therefore for $\xi>0$, we have
  \begin{align*}
    \romnum{1}(\xi)&=-i \int_{\eta>0}\hat{v}(\xi-\eta)\,
    \hat{v}(\eta)\, |\eta|^2d\eta + i
    \int_{\eta<0}\hat{v}(\xi-\eta)\, \hat{v}(\eta)\, |\eta|^2d\eta,\\
    \romnum{2}(\xi)&=-i\, \int_{\eta>0}\hat{v}(\xi-\eta)\,
    \hat{v}(\eta)\, |\eta|^2d\eta \, -\,i\,
    \int_{\eta<0}\hat{v}(\xi-\eta)\, \hat{v}(\eta)\, |\eta|^2d\eta,
  \end{align*}
  and for $\xi<0$,
  \begin{align*}
    \romnum{1}(\xi)&=-i \int_{\eta>0}\hat{v}(\xi-\eta)\,
    \hat{v}(\eta)\, |\eta|^2d\eta + i
    \int_{\eta<0}\hat{v}(\xi-\eta)\, \hat{v}(\eta)\, |\eta|^2d\eta,\\
    \romnum{2}(\xi)&=i\, \int_{\eta>0}\hat{v}(\xi-\eta)\,
    \hat{v}(\eta)\, |\eta|^2d\eta \, +\,i\,
    \int_{\eta<0}\hat{v}(\xi-\eta)\, \hat{v}(\eta)\, |\eta|^2d\eta.
  \end{align*}
  This implies for $\xi>0$,
  \begin{equation*}
    \mathcal{F}(g(v))(\xi)= \romnum{1}(\xi)- \romnum{2}(\xi)=2i\, \int_{\eta<0}\hat{v}(\xi-\eta)\, \hat{v}(\eta)\, |\eta|^2d\eta,
  \end{equation*}
  and for $\xi<0$,
  \begin{equation*}
    \mathcal{F}(g(v))(\xi)= \romnum{1}(\xi)- \romnum{2}(\xi)=-2i\, \int_{\eta>0}\hat{v}(\xi-\eta)\, \hat{v}(\eta)\, |\eta|^2d\eta.
  \end{equation*}
  Using Parseval's relation, we obtain
  \begin{align*}
    \norm{g(v)_x}_2^2 =\norm{\mathcal{F}(g(v)_x)}_2^2 & =
    \int_0^{\infty} |\xi|^2 |\mathcal{F}(g(v))(\xi)|^2 d\xi +
    \int_{-\infty}^{0}   |\xi|^2 |\mathcal{F}(g(v))(\xi)|^2  d\xi\\
    &\le 4 \underbrace{ \int_0^{\infty}   |\xi|^2  \left|  \int_{\eta<0}     \hat{v}(\xi-\eta)  \hat{v}(\eta)  |\eta|^2  d\eta   \right|^2 d\xi  }_{\mathcal{A}}   \\
    &\quad+ 4 \underbrace{ \int_{-\infty}^0 |\xi|^2 \left|
        \int_{\eta>0} \hat{v}(\xi-\eta) \hat{v}(\eta) |\eta|^2 d\eta
      \right|^2 d\xi }_ \mathcal{B}.
  \end{align*}
  Next we estimate $\mathcal{A}$. Note that, for $\eta<0$ and $\xi>0$,
  we have
  \begin{equation*}
    |\eta| \le |\xi-\eta|.
  \end{equation*}
  Using the above inequality we obtain
  \begin{equation*}
    \mathcal{A}\le   \int_0^{\infty}   |\xi|^2  \left(  \int_{\eta<0} \  |\xi-\eta|  |\hat{v}(\xi-\eta)| \  \  |\eta| \ |\hat{v}(\eta)|  d\eta   \right)^2 d\xi.
  \end{equation*}
  Using the Cauchy--Schwarz inequality we find
  \begin{align*}
    \mathcal{A}&\le \int_0^{\infty} |\xi|^2 \left( \int_{\eta<0} \
      |\xi-\eta|^2 |\hat{v}(\xi-\eta)|^2 d\eta \right) \
    \left( \int_{\eta<0} |\eta|^2 \ |\hat{v}(\eta)| ^2 d\eta   \right) d\xi \\
    &\le \norm{v_x}_2^2 \ \ \underbrace{ \int_0^{\infty} |\xi|^2
      \int_{\eta<0} \ |\xi-\eta|^2 |\hat{v}(\xi-\eta)|^2 d\eta
      d\xi}_{\mathcal{C}}.
  \end{align*}
  Next we estimate $\mathcal{C}$. Using change of variables and
  integration by parts, we have
  \begin{align*}
    \mathcal{C}&=\int_0^{\infty}  \xi^2  \left( \int_{\xi}^{\infty}y^2 |\hat{v}(y)| ^2 dy \right)d\xi\\
    &=\frac13 \int_0^{\infty}       \frac{d}{d\xi}(\xi^3)               \left( \int_{\xi}^{\infty}y^2 |\hat{v}(y)| ^2 dy \right)d\xi \\
    &= -  \frac13 \int_0^{\infty}      \xi^3    \           \frac{d}{d\xi} \left( \int_{\xi}^{\infty}y^2 |\hat{v}(y)| ^2 dy \right)d\xi \\
    &= \frac13 \int_0^{\infty} |\xi|^5 |\hat{v}(\xi)|^2 d\xi \le
    \norm{v}_{H^{5/2}}^2.
  \end{align*}
  Therefore we have
  \begin{equation*}
    \mathcal{A}\le \norm{v_x}_2^2\norm{v}_{H^{5/2}}^2. 
  \end{equation*}
  A similar argument shows that
  \begin{equation*}
    \mathcal{B}\le \norm{v_x}_2^2\norm{v}_{H^{5/2}}^2. 
  \end{equation*}
  Finally, we consider
  \begin{align*}
    \norm{H(v_x v_{xx})}_{2}&=  \norm{v_x v_{xx}}_{2}\le   \norm{v_x }_{\infty}\norm{v_{xx}}_{2}\\
    &\le C\norm{v_x }_{H^1}\norm{v}_{H^2}\le
    C\norm{v}_{H^2}\norm{v}_{H^{5/2}}.
  \end{align*}
  This completes the proof of the lemma.
\end{proof}

\begin{proof}[Proof of Theorem \ref{thm:one}]
  The proof compares the error propagation with the exact flow. In our
  approach the necessary regularity for estimating local errors by
  Lemma~\ref{lem:local-1} is ensured by Lemma~\ref{lem:burgers-0}~(i),
  via the following induction argument.

  We make the induction hypothesis that for $k\le n-1$,
  \begin{align*}
    \| u_k \|_{L^2} &\le R, \\
    \| u_k \|_{H^{5/2}} &\le e^{2cRk\dt}\|u_0 \|_{H^{5/2}} \le C_0, \\
    \| u_k - u(t_k) \|_{L^2} &\le \gamma \dt,
  \end{align*}
  where $C_0=e^{2cRT}\|u_0 \|_{H^{5/2}}$ with $c$ from
  Lemma~\ref{lem:burgers-0}~(i), and $\gamma = K(R,T) c_1(C_0) T$ with
  $K(R,T)$ from the local Lipschitz bound \eqref{eq:lip-flow} and
  $c_1(C_0)$ is the constant of Lemma~\ref{lem:local-1} for starting
  values bounded by $C_0$ in $H^{5/2}$.  We then show that the above
  bounds also hold for $k=n$ as long as $n\dt\le T$ and $\dt$ is
  sufficiently small.

  We denote, with $\Phi^t=\Phi^t_{A+B}$ for brevity,
  \begin{equation*}
    u_n^k = \Phi^{(n-k)\dt}(u_k),
  \end{equation*}
  which is the value at time $t_n$ of the exact solution of
  \eqref{eq:main} starting with initial data $u_k$ at time $t_k$. Note
  that
  \begin{equation*}
    u_n = u_n^n, \quad u(t_n)=u_n^0.
  \end{equation*}
  We estimate
  \begin{align*}
    \| u_n - u(t_n) \|_{L^2} &\le \sum_{k=0}^{n-1} \| u_n^{k+1} -
    u_n^k \|_{L^2}
    \\
    &= \sum_{k=0}^{n-1} \| \Phi^{(n-k-1)\dt}(\Pi^\dt(u_k)) -
    \Phi^{(n-k-1)\dt}(\Phi^\dt(u_k)) \|_{L^2}.
  \end{align*}
  For $k\le n-2$, we have $\|\Pi^\dt(u_k)\|_{L^2} = \|u_{k+1}\|_{L^2}
  \le R$, and
  \begin{align*}
    \| \Phi^\dt(u_k) \|_{L^2} &\le \| \Phi^\dt(u_k) - \Phi^\dt(u(t_k))
    \|_{L^2} + \| \Phi^\dt(u(t_k)) \|_{L^2}
    \\
    &\le K(R,\dt)\| u_k - u(t_k) \|_{L^2} + \| u(t_{k+1}) \|_{L^2}
    \\
    &\le K(R,\dt)\gamma \dt + \rho,
  \end{align*}
  cf.~\eqref{eq:HYP}, which is bounded by $R$ if $\dt$ is so small
  that
  \begin{equation*}
    K(R,\dt)\gamma \dt\le R-\rho.
  \end{equation*}
  Using Theorem \ref{theo:main} and Lemma~\ref{lem:local-1} we
  therefore have, for $k\le n-1$ and $n\dt\le T$,
  \begin{align*}
    \| \Phi^{(n-k-1)\dt}(\Pi^\dt(u_k)) &-
    \Phi^{(n-k-1)\dt}(\Phi^\dt(u_k))  \|_{L^2} \\
    &\qquad \le K(R,T) \| \Pi^\dt(u_k)- \Phi^\dt(u_k)\|_{L^2}
    \\
    &\qquad \le K(R,T) c_1(C_0) \dt^2,
  \end{align*}
  where $\Pi$ is the Godunov step operator defined in
  \eqref{eq:godunov}.  With this estimate we obtain, again noting
  $n\dt\le T$,
  \begin{equation*}
    \| u_n - u(t_n) \|_{L^2} \le n K(R,T) c_1(C_0) \dt^2 \le \gamma \dt.
  \end{equation*}
  To prove the boundedness of $u_n$, we choose $\gamma \dt\le
  R-\rho$. Then we have
  \begin{equation*}
    \| u_n \|_{L^2} \le R.
  \end{equation*}
  Since $\|\Phi^t_A(v)\|_{H^{5/2}} \le \| v \|_{H^{5/2}}$, the
  Lemma~\ref{lem:burgers-0}, for $\dt\le \overline t(R)$,
  \begin{equation*}
    \| u_n \|_{H^{5/2}} = \| \Phi^{\dt}_A\circ \Phi^{\dt}_B (u_{n-1})\|_{H^{5/2}}
    \le e^{2cR\dt} \|u_{n-1}\|_{H^{5/2}} \le e^{2cRn\dt} \|u_0\|_{H^{5/2}} .
  \end{equation*}
  Thus, the three necessary results hold by the induction argument,
  and this completes the proof of the theorem.
\end{proof}
\begin{remark} To keep the presentation fairly short we have only
  provided details in the full line case. However, we note that the
  same proofs apply \textit{mutatis mutandis} also in the periodic
  case.
\end{remark}

\section{Strang splitting}
\label{sec:strang}

To prove the correct convergence rate for Strang splitting, we use the
same framework as in the proof of the convergence rate for Godunov
splitting. The major difference between the proofs is that for Strang
splitting we need to use the higher-order midpoint rule, rather than
the rectangle rule applied for the Godunov splitting. In addition, a
higher-order series expansion of the involved terms is also necessary
to obtain the results.

We only present the results in the full line case, and, as before, the
same proofs apply also in the periodic case.  Note that our aim is to
find the error between the operator splitting solution and the exact
(Taylor expanded) solution, and bound it using numerical
quadratures. The proof is longer due to the extra order in the Taylor
expansion.

Also for Strang splitting the proof is based on a local error
estimate.
\begin{lemma}
  \label{lem:local-2}
  Assume that the hypothesis \eqref{eq:HYP1} holds for the solution
  $u(t) = \Phi_{A+B}^t(u_0)$ of \eqref{eq:main}.  If the initial data
  $u_0$ is in $H^{9/2}$, then the local error of the Strang splitting
  \eqref{eq:strang} is bounded in $L^2$ by
  \begin{equation*}
    \norm{\Psi^\dt(u_0)-\Phi_{A+B}^\dt(u_0)}_{L^2} \le c_2 \dt^3,
  \end{equation*}
  where $c_2$ only depends on $\|u_0\|_{H^{9/2}}$.
\end{lemma}
\begin{proof}
  We follow \cite{HoldenLubichRisebro:2011} and use the second-order
  Taylor expansion
  \begin{align*}
    \Phi^\dt_B (v) &= v + \dt B(v) + \tfrac12\dt^2 dB(v)[B(v)]
    \\
    &\quad + \dt^3 \int_0^1 \tfrac12(1-\theta)^2 \Big( d^2B(\Phi^{\theta\dt}_B(v))[B(\Phi^{\theta\dt}_B(v)),B(\Phi^{\theta\dt}_B(v))] \\
    &\qquad\qquad\qquad\qquad\qquad + dB(\Phi^{\theta\dt}_B(v))\big[
    dB(\Phi^{\theta\dt}_B(v))[B(\Phi^{\theta\dt}_B(v))]\big]\Big)\,
    d\theta.
  \end{align*}
  Henceforth we abbreviate the integral remainder term as
  \begin{equation*}
    \dt^3 \int_0^1 \tfrac12(1-\theta)^2 \Big( d^2B(B,B) +dB \, dB \, B
    \Big)\big(\Phi^{\theta\dt}_B(v)\big) \, d\theta .
  \end{equation*}
  Hence,
  \begin{align*}
    u_1&=e^{\dt A}u_0 + \dt e^{\dt A/2}B\big(e^{\dt A/2}u_0\big) +
    \tfrac{1}2\dt^2 e^{\dt A/2}dB\big(e^{\dt A/2}u_0\big)[B\big(e^{\Dt
      A/2}u_0\big)]  \\
    &\quad + \dt^3 \int_0^1 \tfrac12(1-\theta)^2 e^{\dt A/2} \Big(
    d^2B(B,B) +dB \, dB \, B \Big)
    \big(\Phi^{\theta\dt}_B(e^{\dt A/2}u_0)\big) \, d\theta\\
    &= e^{\dt A}u_0 + \dt e^{\dt A/2}B\big(e^{\dt A/2}u_0\big) + e_2,
  \end{align*}
  where $e_2$ is given by
  \begin{align}
    \label{eq:er_new}
    e_2 &:= \tfrac{1}2\dt^2 e^{\dt A/2}dB\big(e^{\dt
      A/2}u_0\big)[B\big(e^{\Dt A/2}u_0\big)] \\
    &\quad + \dt^3 \int_0^1 \tfrac12(1-\theta)^2 e^{\dt A/2} \notag
    \Big( d^2B(B,B) +dB \, dB \, B \Big)
    \big(\Phi^{\theta\dt}_B(e^{\dt A/2}u_0)\big) \, d\theta.
  \end{align}
  Recall \eqref{voc1AA} and \eqref{e1}, viz.
  \begin{equation*}
    u(\dt)=e^{\dt A}u_0 + \int_0^\dt e^{(\dt-s)A} B(e^{sA}u_0)\, ds + e_1
  \end{equation*}
  where
  \begin{equation*}
    e_1 = \int_0^\dt \int_0^s e^{(\dt-s)A}dB(e^{(s-\sigma) A}
    u(\sigma))[e^{(s-\sigma) A}B(u(\sigma))]\, d\sigma \, ds. 
  \end{equation*}
  We express the integrand in $e_1$ by a formula of the type
  \eqref{voc-b} by using
  \begin{equation*}
    G(u(\sigma)) = G(e^{\sigma A}u_0) + \int_0^\sigma
    dG(e^{(\sigma-\tau) A} u(\tau))[e^{(\sigma-\tau) A}B(u(\tau))]\,
    d\tau
  \end{equation*}
  with
  \begin{equation*}
    G(v)=G_{s,\sigma}(v)=dB(e^{(s-\sigma)A}v)[ e^{(s-\sigma)A}\, B(v)],
  \end{equation*}
  and
  \begin{align*}
    dG(v)[w]&=\,d^2B\big(e^{(s-\sigma)A}v\big)
    \big[e^{(s-\sigma)A}w,e^{(s-\sigma)A}B(v)\big] \\ &\quad+
    dB\big(e^{(s-\sigma)A}v\big)\big[e^{(s-\sigma)A} dB(v)[w]\big].
  \end{align*}
  This implies
  \begin{align}
    e_1&=\int_0^{\dt} \int_0^s
    e^{(\dt-s)A}dB\left(e^{sA}u_0\right) \label{eq:er_new1}
    [e^{(s-\sigma)A}B\left(e^{sA}u_0\right)]\,d\sigma ds\\ \notag
    &\quad+ \int_0^{\dt} \int_0^s \int_0^\sigma
    dG_{s,\sigma}\big(e^{(\sigma-\tau)A}u(\tau)\big)\big[e^{(\sigma-\tau)A}
    \notag B(u(\tau))\big] \,d\tau d\sigma ds.
  \end{align}
  We return to the error formula \eqref{err1} and write the principal
  error term
  \begin{align*}
    \dt \, e^{\dt A/2}B\big(e^{\dt A/2}u_0\big) - \int_0^\dt
    e^{(\dt-s)A} B(e^{sA}u_0)\, ds
  \end{align*}
  in second-order Peano form
  \begin{equation*}
    \dt f(\tfrac12\dt) -\int_0^\dt f(s)\,ds = \dt^3 \int_0^1
    \kappa_2(\theta)\, f''(\theta\dt)\, d\theta
  \end{equation*}
  with the second-order Peano kernel $\kappa_2$ of the midpoint rule
  and $f$ is defined by \eqref{f} with
  \begin{equation*}
    f''(s) = e^{(\dt-s)A}[A,[A,B]](e^{sA}u_0).
  \end{equation*} 
  By Lemma~\ref{lem:2COM}, proven below, we obtain the double
  commutator bound
  \begin{align*}
    \norm{[A,[A,B]](v)}_{L^2}\le C\norm{v}_{H^{9/2}}^2.
  \end{align*}
  Thus, it follows that
  \begin{equation*}
    \norm{f''(s)}_{L^2}\le C\norm{u_0}_{H^{9/2}}^2.
  \end{equation*}
  \begin{lemma} \label{lem:2COM} For $v\in H^{9/2}$, we have
    \begin{equation*}
      \norm{[A,[A,B]](v)}_{L^2}\le 2\norm{v_x}_{L^2}\norm{v}_{H^{9/2}}+C\norm{v}_{H^4}^2,
    \end{equation*}
    for some constant $C$.
  \end{lemma}
  \begin{proof}[Proof of Lemma \ref{lem:2COM}]
    The Lie double commutator is given by
    \begin{equation*}
      [A,[A,B]](v):= [A,L](v)= dA(v)[L(v)] - dL(v)[Av],
    \end{equation*}
    where $L$ is defined by
    \begin{align*}
      L(v)&=[A,B](v) = g(v)_x + 2 H(v_x v_{xx}), \\
      \intertext{with} g(v)&:= H(v v_{xx}) - v H(v_{xx}).
    \end{align*}
    A direct computation shows that
    \begin{equation*}
      dL(v)[w]=\Big( H( v w_{xx}+   w v_{xx} )   -       \big( v
      H(w_{xx} )  +  w H(v_{xx}) \big) \Big)_x 
      +2H(v_x w_{xx}+w_xv_{xx}).
    \end{equation*}
    Using the Leibniz rule, we have the following identity
    \begin{align*}
      \big(H(vw)-vH(w)\big)_{xx} &= H(v_{xx}w)+H(w_{xx}v)+2H(v_xw_x)\\
      &\quad -vH(w_{xx})-v_{xx}H(w)-2v_xH(w_x).
    \end{align*}
    Thus
    \begin{equation*}
      dL(v)[w]=\big(H(vw)-vH(w) \big)_{xxx}+\mathcal{E}(v,w),
    \end{equation*}
    where
    \begin{equation*}
      \mathcal{E}(v,w)=\big(v_{xx}H(w)-wH(v_{xx})+2v_xH(w_x)\big)_x.
    \end{equation*}
    For $w=A(v)=H(v_{xx})$, using the property $H^2=-I$, we obtain
    \begin{align*}
      H(vw)-vH(w)&=H(vH(v_{xx})) - v H^2(v_{xx})\\
      &=H(vH(v_{xx})) + v v_{xx}\\
      &=H(vH(v_{xx})-H( v v_{xx}))=-H(g(v)).
    \end{align*}
    Thus,
    \begin{equation*}
      dL(v)[A(v)]=-H\big(g(v)_{xxx} \big)+\mathcal{E}(v,A(v)).
    \end{equation*}
    Again,
    \begin{equation*}
      dA(v)[L(v)] =H\big(L(v)_{xx}\big)=H\big(g(v)_{xxx}\big)-2\big(v_xv_{xx}\big)_{xx}.
    \end{equation*}
    Therefore,
    \begin{equation*} [A,[A,B]](v)=2H\big(g(v)_{xxx}\big)
      +\mathcal{D}(v)
    \end{equation*}
    where
    \begin{equation*}
      \mathcal{D}(v)=-\mathcal{E}(v,A(v))-2\big(v_xv_{xx}\big)_{xx}.
    \end{equation*}
    Repeatedly using that $\norm{v}_{L^\infty}\le \norm{v}_{H^1}$, we
    see that
    \begin{equation*}
      \norm{\mathcal{D}(v)}_2\le C \norm{v}_{H^4}^2,
    \end{equation*}
    for some numerical constant $C$.  Next we claim that
    \begin{equation}\label{d-comm}
      \norm{g(v)_{xxx}}_2\le 2\norm{v_x}_2\norm{v}_{H^{9/2}}.
    \end{equation}
    Using the Parseval relation, we obtain
    \begin{align*}
      \norm{g(v)_{xxx}}_2^2 &=\norm{\mathcal{F}(g(v)_{xxx})}_2^2 \\ &=
      \int_0^{\infty} |\xi|^6 |\mathcal{F}(g(v))(\xi)|^2 d\xi +
      \int_{-\infty}^{0}   |\xi|^6 |\mathcal{F}(g(v))(\xi)|^2  d\xi\\
      &\le 4 \underbrace{ \int_0^{\infty}   |\xi|^6  \left|
          \int_{\eta<0}     \hat{v}(\xi-\eta)  \hat{v}(\eta)  |\eta|^2
          d\eta   \right|^2 d\xi  }_{\mathcal{A}} \\ 
      &\quad + 4 \underbrace{ \int_{-\infty}^0 |\xi|^6 \left|
          \int_{\eta>0} \hat{v}(\xi-\eta) \hat{v}(\eta) |\eta|^2 d\eta
        \right|^2 d\xi }_ \mathcal{B}.
    \end{align*}
    Next we estimate $\mathcal{A}$. Note that, for $\eta<0$ and
    $\xi>0$, we have
    \begin{equation*}
      |\eta| \le |\xi-\eta|.
    \end{equation*}
    Using the above inequality we obtain
    \begin{equation*}
      \mathcal{A}\le   \int_0^{\infty}   |\xi|^6  \left(
        \int_{\eta<0} \  |\xi-\eta|  |\hat{v}(\xi-\eta)| \  \  |\eta|
        \ |\hat{v}(\eta)|  d\eta   \right)^2 d\xi. 
    \end{equation*}
    Using the Cauchy--Schwarz inequality we infer
    \begin{align*}
      \mathcal{A}&\le \int_0^{\infty} |\xi|^6 \left( \int_{\eta<0} \
        |\xi-\eta|^2 |\hat{v}(\xi-\eta)|^2 d\eta \right) \
      \left( \int_{\eta<0} |\eta|^2 \ |\hat{v}(\eta)| ^2 d\eta   \right) d\xi \\
      &\le \norm{v_x}_2^2 \ \ \underbrace{ \int_0^{\infty} |\xi|^6
        \int_{\eta<0} \ |\xi-\eta|^2 |\hat{v}(\xi-\eta)|^2 d\eta
        d\xi}_{\mathcal{C}}.
    \end{align*}
    Next we estimate $\mathcal{C}$. Using a change of variables and
    integration by parts, we have
    \begin{align*}
      \mathcal{C}&=\int_0^{\infty}  \xi^6  \left(
        \int_{\xi}^{\infty}y^2 |\hat{v}(y)| ^2 dy \right)d\xi\\ 
      &=\frac17 \int_0^{\infty}       \frac{d}{d\xi}(\xi^7)
      \left( \int_{\xi}^{\infty}y^2 |\hat{v}(y)| ^2 dy \right)d\xi \\ 
      &= -  \frac17 \int_0^{\infty}      \xi^7    \
      \frac{d}{d\xi} \left( \int_{\xi}^{\infty}y^2 |\hat{v}(y)| ^2 dy
      \right)d\xi \\ 
      &= \frac17 \int_0^{\infty} |\xi|^9 |\hat{v}(\xi)|^2 d\xi \le
      \norm{v}_{H^{9/2}}^2.
    \end{align*}
    Therefore we have
    \begin{equation*}
      4 \mathcal{A}\le \norm{v_x}_2^2\norm{v}_{H^{9/2}}^2. 
    \end{equation*}
    A similar argument shows that
    \begin{equation*}
      4 \mathcal{B}\le \norm{v_x}_2^2\norm{v}_{H^{9/2}}^2. 
    \end{equation*}
    This completes the proof of \eqref{d-comm} and thereby of the
    lemma.
  \end{proof}
  Now for the difference of \eqref{eq:er_new} and \eqref{eq:er_new1},
  \begin{equation*}
    e_2 -e_1 = \tfrac12 \dt^2g(\tfrac12\dt,\tfrac12\dt) -
    \int_0^\dt\int_0^s g(s,\sigma)\,d\sigma\, ds + \tilde e_2 - \tilde
    e_1,
  \end{equation*}
  where
  \begin{align*}
    g(s,\sigma) &= e^{(\dt-s)A} dB( e^{s A}u_0)\,
    [e^{(s-\sigma)A}B(e^{\sigma A}u_0)], \\
    \tilde{e}_1 &= \int_0^{\dt} \int_0^s \int_0^\sigma
    dG_{s,\sigma}\left(e^{(\sigma-\tau)A}u(\tau)\right)
    e^{(\sigma-\tau)A} B(u(\tau)) \,d\tau d\sigma ds, \\
    \tilde{e}_2 &= \dt^3 \int_0^1 \tfrac12(1-\theta)^2 e^{\dt A/2}
    \Bigl( d^2B(B,B) + dB \,dB \, B\Bigr)(\Phi^{\theta\dt}_B(u_0))\,
    d\theta.
  \end{align*}
  To estimate the remainder terms $\tilde e_i$, for $i=1,2$, we
  calculate
  \begin{align*}
    \norm{dG_{s,\sigma}(v)w}_{L^2}&\le
    \norm{d^2B\left(e^{(s-\sigma)A}v\right)[e^{(s-\sigma)A}B(v),e^{(s-\sigma)A}
      w]}_{L^2}\\
    &\quad+ \norm{dB\left(e^{(s-\sigma)A} v
      \right)[e^{(s-\sigma)A}dB(v)[w]]}_{L^2}\\
    &\le C\left(\norm{B(v)}_{H^{1}}\norm{w}_{H^{1}}+\norm{v}_{H^{1}}
      \norm{dB(v)[w]}_{H^{1}}\right)\\
    &\le C\left(\norm{v}_{H^{2}}^2 \norm{w}_{H^{1}} +
      \norm{v}_{H^{1}}\norm{v}_{H^{2}} \norm{w}_{H^{2}}\right)\\
    &\le C\norm{v}_{H^{2}}^2\norm{w}_{H^{2}},
  \end{align*}
  and
  \begin{align*}
    \norm{\left(d^2B(B,B) + dB\, dB\, B\right)(v)}_{L^2}&\le
    \norm{d^2B(B(v),B(v))}_{L^2}+ \norm{dB(v)[dB(v)[B(v)]]}_{L^2}\\
    &\le C\left(\norm{B(v)}_{H^{1}}^2 +
      \norm{v}_{H^{1}}\norm{dB(v)[B(v)]}_{H^{1}}\right) \\
    &\le C\left(\norm{v}_{H^{2}}^4 +
      \norm{v}_{H^{2}}^2\norm{B(v)}_{H^{2}}\right) \\
    &\le C \norm{v}_{H^{3}}^4.
  \end{align*}
  Then, using Lemma 3.1, the remainder terms are bounded by
  \begin{equation}
    \label{fbound1}
    \norm{ \tilde e_1}_{L^2} + \norm{\tilde e_2}_{L^2} \le C \dt^3
    \norm{u_0}_{H^3}^4.
  \end{equation}
  The first two terms in $e_2-e_1$ are the quadrature error of a
  first-order two-dimensional quadrature formula, which is bounded by
  \begin{multline*}
    \Bigl\| \tfrac12 \dt^2g(\tfrac12\dt,\tfrac12\dt) -
    \int_0^\dt\int_0^s g(s,\sigma)\,d\sigma\, ds \Bigr\|_{L^2}
    \\
    \le C \dt^3 \bigl( \max \| \partial g/\partial s\|_{L^2} + \max
    \| \partial g/\partial\sigma \|_{L^2}\bigr),
  \end{multline*}
  where the maxima are taken over the triangle
  $\{(s,\sigma):0\le\sigma\le s \le \dt\}$.  In order to estimate the
  partial derivatives we write
  \begin{equation*}
    g(s,\sigma)=e^{(\dt-s)A} dB(v(s)) w(s,\sigma),
  \end{equation*}
  where
  \begin{equation*}
    v(s)=e^{sA}u_0\ \text{ and }\ w(s,\sigma)=e^{(s-\sigma)A}B(v(\sigma)).
  \end{equation*}
  With this notation
  \begin{align*}
    \frac{\partial g}{\partial s}&= e^{(\dt-s)A}
    \left(-AdB(v(s))w(s,\sigma)) +d^2B(Av(s),w(s,\sigma))+
      dB(v(s))Aw(s,\sigma)\right)\\
    &= e^{(\dt-s)A}\left(-A(v(s)w(s,\sigma))+Av(s) w(s,\sigma)+
      v(s)Aw(s,\sigma)\right)_x.
  \end{align*}
  Now
  \begin{align*}
    \norm{\frac{\partial g}{\partial
        s}}_{L^2}&\le\norm{-A(v(s)w(s,\sigma)) + (Av(s))w(s,\sigma) +
      v(s)Aw(s,\sigma)}_{H^1}\\
    &\le \norm{A(vw)}_{H^1}+ \norm{(Av)w}_{H^1} + \norm{vAw}_{H^1}\\
    &\le C\norm{v}_{H^3}\norm{w}_{H^3}\le C \norm{v}^3_{H^4}\le C
    \norm{u_0}^3_{H^4}.
  \end{align*}
  For the other partial derivative, we get
  \begin{align*}
    \frac{\partial g}{\partial \sigma} &= e^{(\dt-s)A}
    dB(v(s))\left(e^{(s-\sigma)A}\left(-AB(v(s))
        +dB(v(\sigma))Av(\sigma)\right)\right),
  \end{align*}
  so that
  \begin{align*}
    \norm{ \frac{\partial g}{\partial \sigma} }_{L^2}&\le C
    \norm{v(s)}_{H^{1}} \norm{-AB(v(s))
      +dB(v(\sigma))Av(\sigma)}_{H^{1}}\\
    &\le  \norm{v}_{H^1}\norm{AB(v)}_{H^1}+\norm{v}_{H^1}\norm{dB(v)[Av]}_{H^1}\\
    &\le \norm{v_{H^1}}\left(\norm{{(vv_x)}_{xx}}_{H^1}
      +\norm{v}_{H^2}\norm{A(v)}_{H^2}\right)\\
    &\le C \norm{v}_{H^4}^3.
  \end{align*}
  Therefore
  \begin{equation*}
    \norm{ \frac{\partial g}{\partial \sigma} }_{L^2}\le C \norm{u_0}_{H^4}^3.
  \end{equation*}
  Thus
  \begin{equation}
    \label{fbound3}
    \norm{e_2 - e_1}_{L^2} \le C \dt^3\big( \norm{u_0}_{H^4}^3+
    \norm{u_0}_{H^3}^4\big) \le C \dt^3, 
  \end{equation}
  which together with the bound for the quadrature error of the
  midpoint rule for $f$ yields the stated result.
\end{proof}
\begin{proof}[Proof of Theorem~\ref{thm:two}] 
  We argue as in the proof of Theorem~\ref{thm:one}, but now assume
  inductively that $\norm{u_k-u(t_k)}_{L^2}\le \gamma\Dt^2$. With
  $\Psi$ denoting the Strang step operator defined in
  \eqref{eq:strang}, we have
  \begin{align*}
    \norm{u_n - u(t_n)}_{L^2} &\le \sum_{k=0}^{n-1}
    \norm{\Phi^{(n-k-1)\dt}(\Psi^\dt(u_k)) -
      \Phi^{(n-k-1)\dt}(\Phi^\dt(u_k))}_{L^2}
    \\
    &\le \sum_{k=0}^{n-1} K(R,T) \| \Psi^\dt(u_k)-\Phi^\dt(u_k)
    \|_{L^{2}}
    \\
    &\le \sum_{k=0}^{n-1} K(R,T) c_2(C_0) \dt^3 \le K(R,T) c_2(C_0) T
    \dt^2.
  \end{align*}
  This completes the proof of Theorem~\ref{thm:two}.
\end{proof}



\bigskip \textbf{Acknowledgement.}  HH is grateful to Luc Molinet for
helpful discussions.



\begin{thebibliography}{99}
\bibitem{abdel} L. ~Abdelouhab, J. ~L. Bona, M. ~Felland, and
  J.~C. Saut.  \newblock Nonlocal models for nonlinear dispersive
  waves.  \newblock \textit{Physica D 40} 360--392 (1989).

\bibitem{AmbrosettiProdi} A. Ambrosetti and G. Prodi.  \newblock
  \textit{A Primer of Nonlinear Analysis}.  \newblock Cambridge UP,
  Cambridge (1995).
%
%
\bibitem{ablowitz} M. ~J. Ablowitz and A. ~S. Fokas.  \newblock The
  inverse scattering transform for the Benjamin--Ono equation, a pivot
  for multidimensional problems.  \newblock \textit{Stud. Appl. Math.}
  68:1--10 (1983).
%
%
\bibitem{benjamin} T.~B. Benjamin.  \newblock Internal waves of
  permanent form in fluid of great depth.  \newblock
  \textit{J. Fluid. Mech.}  29:559--592 (1967).
%
%
%

\bibitem{burq} N.~Burq and F.~Planchon.  \newblock On well-posedness
  for the Benjamin--Ono equation.  \newblock \textit{Math. Ann.}
  340:497--542 (2008).
%
\bibitem{HoldenKarlsenKarper:2011} H. Holden, K. H. Karlsen,
  T. Karper.  \newblock Operator splitting for two-dimensional
  incompressible fluid equations.  \newblock \textit{Math. Comp.}
  82:719--748 (2013).
%
\bibitem{HoldenKarlsenKarper:2013} H. Holden, K. H. Karlsen,
  T. Karper.  \newblock Operator splitting for well-posed active
  scalar equations.  \newblock \textit{SIAM J. Math. Anal.}
  45(1):152--180 (2013).
%
\bibitem{Holden:book} H. Holden, K.~H. Karlsen, K.-A. Lie, and
  N.~H. Risebro.  \newblock {\em Splitting for Partial Differential
    Equations with Rough Solutions. Analysis and Matlab Programs.}
  \newblock European Math. Soc. Publishing House, Z\"urich (2010).
%
\bibitem{HoldenKarlsenRisebro:1999} H. Holden, K. H. Karlsen, and
  N. H. Risebro.  \newblock Operator splitting methods for generalized
  {K}orteweg--de {V}ries equations.  \newblock
  \textit{J. Comput. Phys.} 153:203--222 (1999).
%
\bibitem{HoldenKarlsenRisebroTao:2009} H. Holden, K. H. Karlsen,
  N. H. Risebro, and T. Tao.  \newblock Operator splitting methods for
  the {K}orteweg--de {V}ries equation.  \newblock \textit{Math. Comp.}
  80:821--846 (2011).
 %
%
\bibitem{HoldenLubichRisebro:2011} H. Holden, C. Lubich, and
  N. H. Risebro.  \newblock Operator splitting for partial
  differential equations with Burgers nonlinearity.  \newblock
  \textit{Math. Comp.} 82:173--185 (2013).
%
%
%
\bibitem{ionescu} A. ~Ionesc and C. ~E. Kenig.  \newblock Global
  well posedness of the Benjamin--Ono equation in low regularity
  spaces.  \newblock \textit{J. Amer. Math. Soc.}, 20:753--798 (2007).


\bibitem{iorio} R. ~Iorio.  \newblock On the Cauchy problem for the
  Benjamin--Ono equation.  \newblock \textit{Comm. Part. Diff. Eq.},
  11:1031--1081 (1986).

%


\bibitem{KochTzvetkov2005} H. Koch, and N. Tzvetkov.  \newblock
  Nonlinear wave interactions for the Benjamin--Ono equation.
  \newblock \textit{International Mathematics Research Notices.} No-30
  (2005).

\bibitem{Leela} V. Lakshmikantham and S. Leela.  \newblock {\em
    Nonlinear Differential Equations in Abstract Spaces.}  \newblock
  Pergamon Press, Oxford (1981).
 


%
%
\bibitem{Lubich:2008} C. Lubich.  \newblock On splitting methods for
  Schr\"odinger--Poisson and cubic nonlinear Schr\"odinger equations.
  \newblock \textit{Math. Comp.}  77:2141--2153 (2008).
  
%

\bibitem{Mol1} L.~Molinet.  \newblock Global well-posedness in
$L^2$ for the periodic Benjamin--Ono equation.  
\newblock \textit{American Journal of Mathematics.} Vol. 130, No. 3: 635-683 (2008).


\bibitem{molinet3} L.~Molinet.  \newblock Global well-posedness in the
  energy space for the Benjamin--Ono equation on the circle.
  \newblock \textit{Math. Ann.}  337:353--383 (2007).

\bibitem{ono} H.~Ono.  \newblock Algebraic solitary waves in
  stratified fluids.  \newblock \textit{J. Phy. Soc. Japan}
  39(4):1082--1091 (1975).



\bibitem{ponce} G.~Ponce.  \newblock On the global well posedness of
  the Benjamin--Ono equation.  \newblock \textit{Diff. Int. Eq.}
  4:527--542 1991).










%



\bibitem{Tao:2004} T. Tao.  \newblock Global well-posedness of the
  Benjamin--Ono equation in $H^1(\R)$.  \newblock \textit{
    J. Hyp. Diff. Equations} 1(1):27--49 (2004).

%





\end{thebibliography}
\end{document}